\documentclass[12pt]{amsart}
\usepackage{ amsmath, amsthm, amsfonts, amssymb, color, xcolor}
 \usepackage{mathrsfs}
\usepackage{amsfonts, amsmath}
 \usepackage{amsmath,amstext,amsthm,amssymb,amsxtra,verbatim}

\begin{document}

\newfam\msbfam
\font\tenmsb=msbm10    \textfont\msbfam=\tenmsb \font\sevenmsb=msbm7
\scriptfont\msbfam=\sevenmsb \font\fivemsb=msbm5
\scriptscriptfont\msbfam=\fivemsb
\def\Bbb{\fam\msbfam \tenmsb}
\newfam\bigfam
\font\tenbig=msbm10 scaled \magstep2   \textfont\bigfam=\tenbig
\font\sevenbig=msbm7 scaled \magstep2 \scriptfont\bigfam=\sevenbig
\font\fivebig=msbm5 scaled \magstep2

\scriptscriptfont\bigfam=\fivebig

\def\dsum{\displaystyle\sum}
\def\dfrac{\displaystyle\frac}
\def\dlim{\displaystyle\lim}
\def\dint{\displaystyle\int}
\def\dsup{\displaystyle\sup}
\def\nln{\newline}

\newtheorem{thm}{Theorem}[section]
\newtheorem{lem}{Lemma}[section]
\newtheorem{prop}{Proposition}[section]
\newtheorem{rem}{Remark}[section]
\newtheorem{cor}{Corollary}[section]
\newtheorem{defn}{Definition}[section]
\newtheorem{pf}{\it Proof.}
\renewcommand\thepf{}

\newcommand{\beq}{\begin{equation}}
\newcommand{\eeq}{\end{equation}}
\numberwithin{equation}{section}

\title[ Sharp inequalities for  martingales]
{Sharp inequalities for linear combinations of orthogonal martingales}

\author{Yong Ding}

\address{Yong Ding, Laboratory of Mathematics and Complex Systems, School of Mathematical Sciences, Beijing Normal University, Ministry of Education of China, Beijing 100875, China }
\email{dingy@bnu.edu.cn }

\author{Loukas Grafakos}

\address{Loukas Grafakos, Department of Mathematics, University of Missouri, Columbia MO 65211, USA}
\email{grafakosl@missouri.edu}

\author{Kai Zhu${}^1$}

\address{Kai Zhu, School of Mathematical Sciences, Beijing Normal University, Beijing 100875, China}
\email{kaizhu0116@126.com }

\thanks{{\it Keywords and phases:} Sharp inequalities, Martingales, orthogonal, Hilbert transform, discrete.}
\thanks{The first author is supported by NSFC(11371057,11471033,11571160)  and the Fundamental Research Funds for the Central Universities(2014KJJCA10).}
\thanks{The second author  would like to acknowledge the support of Simons Foundation and of the
University of Missouri Research Board and Research Council.}
\thanks{${}^1$Corresponding author.}

\begin{abstract}
For any two real-valued continuous-path martingales $X=\{X_t\}_{t\geq 0}$ and $Y=\{Y_t\}_{t\geq 0}$,  with $X$ and $Y$ being orthogonal and $Y$ being differentially subordinate to $X$, we obtain sharp $L^p$ inequalities for martingales of the form $aX+bY$ with $a, b$ real numbers. The best $L^p$ constant is equal to the norm of the operator $aI+bH$ from $L^p$ to $L^p$, where $H$ is the Hilbert transform on the circle or real line. The   values of these norms  were found by Hollenbeck, Kalton and Verbitsky \cite{HKV}.
\end{abstract}

\maketitle

\section{Introduction}
The research on martingale inequalities was initiated in 1966 by Burkholder \cite{B}  and was
further pursued  in \cite{B0}, \cite{B1} and \cite{B2}, where techniques for  sharp estimates for them    were developed. Martingale inequalities nowdays find applications in probability and analysis and their impact is quite far-reaching.

Based on the techniques of Burkholder,
 Ba\~{n}uelos and Wang \cite{BW} obtained   sharp inequalities for orthogonal
 martingales, and used them to
provide  probabilistic proofs to the results of  Pichorides \cite{Pi}  concerning  the
norm of the Hilbert transform on $L^p(\mathbb R)$ and of Iwaniec and Martin \cite{IM}
about the norm of the   Riesz transforms on $L^p(\mathbb R^n)$, $1<p<\infty$.

We  describe the pertinent framework for this paper.
Let $(\Omega, \mathcal{F}, P)$ be a probability space and $\mathcal{F}=\{\mathcal{F}_t\}_{t\geq 0}$ be a nondecreasing family of sub-$\sigma$-fields of
$\mathcal{F}_\infty  = \cup_{t\ge 0} \mathcal{F}_t $. Let $X=\{X_t\}_{t\geq 0}$ and $Y=\{Y_t\}_{t\geq 0}$ be two real-valued martingales with respect to $\mathcal{F}$. We say that $X$ is orthogonal to $Y$ if $\langle X, Y \rangle_t=0$ for all $t\geq 0$, where $\langle X, Y \rangle_t$ is the quadratic covariation between $X$ and $Y$. We also say that $Y$ is differentially subordinate to $X$ (see \cite{BW}) if $\langle X\rangle_t-\langle Y\rangle_t$ is a nondecreasing function of $t$ for $t\geq 0$, where $\langle X\rangle_t$ is the quadratic variation of $X$.

For $1<p<\infty$, define $n_p=\cot (\pi/(2p^\ast))$, where $p^\ast=\max (p,p/(p-1))$.
This constant is exactly the operator norm of Hilbert transform $H$ on $L^p(\mathbb R)$
and of the conjugate function  $H^{\mathbb T}$ on $L^p(\mathbb T)$,
where $T$ is the unit circle; see  Pichorides \cite{Pi}.

For continuous-path real-valued martingale $X=\{X_t\}_{t\geq 0}$, $1<p<\infty$, define
$$
\|X\|_p=\sup_{t\geq 0}\|X_t\|_p,
$$
where $\|X_t\|_p=(E|X_t|^p)^{1/p}$. In \cite{BW}, Ba\~{n}uelos and Wang obtained the following   result:

\bigskip
\noindent \textbf{Theorem A.} (\cite{BW}) \emph{Let $X$ and $Y$ be two real-valued continuous-path martingales such that $X$ and $Y$ are orthogonal and $Y$ is differentially subordinate to $X$. Then for $1<p<\infty$,
\beq\label{(BW1)}
\|Y\|_p\leq n_p\|X\|_p
\eeq
and
\beq\label{(BW2)}
\|(X^2+Y^2)^{1/2}\|_p\leq E_p\|X\|_p,
\eeq
where $E_p=(1+n_p^2)^{1/2}$. The constants are best possible.}

The results of Theorem A are the martingale analogues of the results of Pichorides \cite{Pi},
Iwaniec and Martin \cite{IM}, and Ess\'{e}n \cite{E}.

\smallskip

For $a, b \in \mathbb R$, $1<p<\infty$, define
\beq\label{(0.1)}
B_p=\max_{x\in \Bbb R}\frac{|ax-b+(bx+a)\tan \gamma|^p+|ax-b-(bx+a)\tan \gamma|^p}{|x+\tan \gamma|^p+|x-\tan \gamma|^p},
\eeq
where $\gamma =\frac{\pi}{2p}$. $B_p$ can be  equivalently defined as
\beq\label{(1.00)}
B_p=(a^2+b^2)^{p/2}\max_{0\leq \theta \leq 2\pi}\frac{|\cos(\theta +\theta_0)|^p+|\cos(\theta +\theta_0+\frac{\pi}{p})|^p}{|\cos \theta|^p+|\cos(\theta +\frac{\pi}{p})|^p},
\eeq
and
\beq\label{(1.0)}
B_p=(a^2+b^2)^{p/2}\max_{0\leq \vartheta \leq 2\pi}\frac{|\cos(\vartheta -\theta_0)|^p+|\cos(\vartheta -\theta_0+\frac{\pi}{p})|^p}{|\cos \vartheta|^p+|\cos(\vartheta +\frac{\pi}{p})|^p},
\eeq
where $\tan \theta_0=b/a$. These constants appeared in the work of
Hollenbeck, Kalton and Verbitsky \cite{HKV} who showed that the norm of $aI+bH^{\mathbb T}$ from $L^p(\mathbb T)$ to $L^p(\mathbb T)$ is equal to $B_p^{1/p}$, where $I$ is the identity operator and $H^{\mathbb T}$ is the conjugate function operator on the circle. The same assertion is also true for the norm of $aI+bH$ from $L^p(\mathbb R)$ to $L^p(\mathbb R)$, where $H$ is the Hilbert transform on real line, through a   dilation argument known as ``blowing up the circle" (see \cite{Z}, Chapter XVI, Theorem 3.8).  Recently, Ding, Grafakos and Zhu \cite{DGZ} provided a direct proof  of the sharp $L^p(\mathbb R)$ inequality for $aI+bH$  by
an argument that uses  an explicit formula for a crucial  subharmonic majorant.

In this work, we prove sharp inequalities for the martingale  $aX+bY$, where $X$ and $Y$ are as in
Theorem A and  $a, b$ are   arbitrary real numbers.
Motivated by the usefulness of the explicit formula
of the crucial subharmonic majorant $G$ in \cite{DGZ}, we derive two alternative
explicit expressions for this function  (Lemma~\ref{Lemma 2}),
and use them appropriately   in the proof of the main estimate \eqref{(1.1)} below.

\begin{thm}\label{T1.1} Let $X$ and $Y$ be two real-valued continuous-path martingales such that $X$ and $Y$ are orthogonal and $Y$ is differentially subordinate to $X$. Let $B_p$ be given by \eqref{(1.0)}. Then for $a, b\in \mathbb R$ and  $1<p<\infty$ we have
\beq\label{(1.1)}
\|aX+bY\|_p\leq B_p^{1/p}\|X\|_p.
\eeq
The constant  $B_p^{1/p}$ is the best possible in this inequality.
\end{thm}

Inequality \eqref{(1.1)} is the martingale analogue of that in Hollenbeck et al. \cite{HKV} for analytic functions in the unit disc.

We now turn to the proof of this theorem. Without loss of generality, we assume that $a=\cos \theta_0, b=\sin \theta_0$, so that $a^2+b^2=1$. We also assume throughout the paper that $X_0=Y_0=0$.

\section{Some Lemmas}
In this section we discuss some crucial lemmas  in the proof of the main theorem. The first lemma is a version of Lemma 4.2 in \cite{HKV},   in which we derive an explicit formula for a  subharmonic  function $G$ that plays a crucial role in the proof.

\begin{lem}\label{Lemma 1}\cite[Lemma 3.2]{DGZ} Let $1<p<\infty$, $B_p$ be given by \eqref{(1.0)}, $T=\{re^{it}:r>0,t_0<t<t_0+\frac{\pi}{p}\}$, where $t_0$ is the value that makes right part of \eqref{(1.0)} attain its maximum, and there exists $\varepsilon>0$ such that $t_0-\varepsilon<t_0<t_0+\pi/p<t_0+\pi-\varepsilon$. Let $z=re^{it}, z_0=re^{it_0}, G(z)=G(re^{it})$ be $\pi$-periodic of $t$ and when $t_0-\varepsilon<t<t_0+\pi-\varepsilon$:
$$
G(z)=\begin{cases}
B_p|\mathrm{Re}z_0|^{p-1}\textup{sgn}(\mathrm{Re}z_0)\mathrm{Re}[(\frac{z}{z_0})^pz_0]-|a\mathrm{Re}z_0+b\mathrm{Im}z_0|^{p-1}\\
\quad \quad\times\,\, \textup{sgn}(a\mathrm{Re}z_0+b\mathrm{Im}z_0)(a\mathrm{Re}[(\frac{z}{z_0})^pz_0]+b\mathrm{Im}[(\frac{z}{z_0})^pz_0]), &\textup{if}\quad z\in T  \\
B_p|\mathrm{Re}z|^p-|a\mathrm{Re}z+b\mathrm{Im}z|^p, \quad \quad \quad\quad\quad\quad\quad\quad \,\,\qquad &\textup{if}\quad z \notin T.
\end{cases}
$$
Then $G(z)$ is subharmonic on $\Bbb C$ and satisfies
\beq\label{(5)}
|a\, \mathrm{Re}z+b\, \mathrm{Im}z|^p\leq B_p|\mathrm{Re}z|^p-G(z).
\eeq
for all $z\in \Bbb C$.
\end{lem}

In the next lemma, we  provide two other explicit formulas for $G$ centered around the points  $t_0$ and $u_0=t_0+\pi/p$, respectively.

\begin{lem}\label{Lemma 2} Let $1<p<\infty$, $B_p, T$ and $t_0, \varepsilon$ be as in Lemma~\ref{Lemma 1}. Let $z=re^{it}, z_0=re^{it_0}$. Then for $z=re^{it}\in T$, $G(z) $ in Lemma~\ref{Lemma 1} has the following equivalent expressions:
\beq\label{(5.1)}
G(z)=
r^p\bigg[B_p \dfrac{|\cos t_0|^p}{\cos t_0}\cos(p(t-t_0)+t_0)-\dfrac{|\cos (t_0-\theta_0)|^p}{\cos (t_0-\theta_0)}\cos(p(t-t_0)+t_0-\theta_0)\bigg]
\eeq
and
\beq\label{(5.2)}
G(z)=
r^p\bigg[B_p \dfrac{|\cos u_0|^p}{\cos u_0}\cos(p(t\!-\! u_0)\!+\! u_0)-\dfrac{|\cos (u_0\!-\! \theta_0)|^p}{\cos (u_0\! -\! \theta_0)}\cos(p(t\!- \! u_0)\! +\! u_0\! -\! \theta_0)\bigg],
\eeq
where $u_0=t_0+\pi/p$, $\tan \theta_0=b/a$, $G(z)$ is $\pi$-periodic of $t$ and $t_0-\varepsilon <t<t_0+\pi-\varepsilon$.
\end{lem}

\emph{Proof.} Expression \eqref{(5.1)} is just the one  given in Lemma 3.2 in \cite{DGZ}.
We now  prove \eqref{(5.2)}. In the proof of Lemma 3.2 in \cite{DGZ}, using the notation in that
reference, we have
\beq\label{(5.1.1)}
h(x)=\widetilde{f}(\widetilde{p}\widetilde{t_0})\cos (x-\widetilde{p}\widetilde{t_0})+\widetilde{f}'_+(\widetilde{p}\widetilde{t_0})\sin (x-\widetilde{p}\widetilde{t_0}),
\eeq
where $\widetilde{p}=p/2, \widetilde{t_0}=2t_0$ and
\beq\label{(5.1.1.1)}
\widetilde{f}(t)=B_p|\cos(t/p)|^p-|a\cos(t/p)+b\sin(t/p)|^p,
\eeq
if we can prove
\beq\label{(5.1.2)}
h(x)=\widetilde{f}(\widetilde{p}\widetilde{t_0}+\pi)\cos (x-\widetilde{p}\widetilde{t_0}-\pi)+\widetilde{f}'_+(\widetilde{p}\widetilde{t_0}+\pi)\sin (x-\widetilde{p}\widetilde{t_0}-\pi),
\eeq
then following the proof of Lemma 3.2 in \cite{DGZ}, we deduce \eqref{(5.2)} when $z\in T $.

To obtain  \eqref{(5.1.2)}, in view of \eqref{(5.1.1)},  it is sufficient to show that
\beq\label{(5.1.3)}
\widetilde{f}(\widetilde{p}\widetilde{t_0})+\widetilde{f}(\widetilde{p}\widetilde{t_0}+\pi)=0
\eeq
and
\beq\label{(5.1.4)}
\widetilde{f}'_+(\widetilde{p}\widetilde{t_0})+\widetilde{f}'_+(\widetilde{p}\widetilde{t_0}+\pi)=0.
\eeq
In fact,
$$
\widetilde{f}(\widetilde{p}\widetilde{t_0}+\pi)=\widetilde{f}(pt_0+\pi)=B_p|\cos(t_0+\pi/p)|^p-
|a\cos(t_0+\pi/p)+b\sin(t_0+\pi/p)|^p,
$$
by
\beq\label{(5.1.5)}
B_p=(a^2+b^2)^{p/2}\frac{|\cos(t_0 -\theta_0)|^p+|\cos(t_0 -\theta_0+\frac{\pi}{p})|^p}{|\cos t_0|^p+|\cos(t_0 +\frac{\pi}{p})|^p},
\eeq
where $\tan \theta_0=b/a$, we have
$$
\begin{array}{cl}
&B_p|\cos(t_0+\pi/p)|^p-|a\cos(t_0+\pi/p)+b\sin(t_0+\pi/p)|^p\\
=&-B_p|\cos t_0|^p+|a\cos t_0+b\sin t_0|^p
=-\widetilde{f}(pt_0)=-\widetilde{f}(\widetilde{p}\widetilde{t_0}),
\end{array}
$$
so we get \eqref{(5.1.3)}.

For \eqref{(5.1.4)}, note that for $1<p<\infty$, $\widetilde{f}(t)$ is $p\pi$-periodic and continuously differentiable. By \eqref{(1.0)}, $g(t)=\widetilde{f}(t)+\widetilde{f}(t+\pi)\geq0$ and $g(t)$ has a minimum at $\widetilde{p}\widetilde{t_0}$, so
$$
\widetilde{f}'_+(\widetilde{p}\widetilde{t_0})+\widetilde{f}'_+(\widetilde{p}\widetilde{t_0}+\pi)=
g'(\widetilde{p}\widetilde{t_0})=0,
$$
Thus the lemma is proved. \qed

The preceding lemma indicates that the   function $G$ has some symmetry properties in terms of $t_0$ and $u_0$.

\section{Proof of Theorem~\ref{T1.1}}\label{Proof}
The proof of Theorem~\ref{T1.1} is based on the techniques of Burkholder; also see \cite{BW}. We choose the appropriate function for Theorem~\ref{T1.1} to be the opposite   of function $G$ in Lemma~\ref{Lemma 1} and use the explicit formulas for $G$ obtained by Lemma~\ref{Lemma 2}.

For $x,y\in \mathbb R$, $1<p<\infty$, set
$$
V(x,y)=|ax+by|^p-B_p|x|^p,
$$
where $x=r\cos t, y=r\sin t$, and $0<t<2\pi$. Define
$$
U(x,y)=-G(x+iy)=-G(z),
$$
where $z=re^{it}$ and $G(z)$ is the function in Lemma~\ref{Lemma 1}. Then by Lemma~\ref{Lemma 1}, we have
\beq\label{(3.1)}
V\leq U.
\eeq

Denoting by $U_{xx}, U_{yy}$   the second order partial derivatives of $U(x,y)$, we need only to show that for all $h, k \in \mathbb R$,
\beq\label{(3.2)}
U_{xx}(x,y)h^2+U_{yy}(x,y)k^2\leq -c(x,y)(h^2-k^2)
\eeq
for $(x,y)\in S_i$, where $S_i, i\geq 1$ is a sequence of open connected sets such that the union of the closure of $S_i$ is $\mathbb R^2$, and $c(x,y)\geq 0$ that is bounded on $1/\delta \leq r\leq \delta$ for any $\delta>0$. In fact, using Proposition 1.2 with Remark 1.1 in \cite{BW}, by \eqref{(3.2)}, we can get
$$
EV(X_t,Y_t)\leq EU(X_t,Y_t)\leq EU(X_0,Y_0)\leq 0,
$$
thus
$$
E|aX_t+bY_t|^p\leq B_pE|X_t|^p,
$$
then we get \eqref{(1.1)}.

To show \eqref{(3.2)},we split the argument into two cases. First, for $z=x+iy \notin T $ we have
$$
U(x,y)=|ax+by|^p-B_p|x|^p,
$$
and by a direct calculation we obtain from this that
\beq
U_{xx}(x,y)=p(p-1)\bigg(|ax+by|^{p-2}a^2-B_p|x|^{p-2} \bigg)
\eeq
except on the lines $\{z:\,x=0\}$ and $\{z:\,ax+by=0\}$,   and
\beq
U_{yy}(x,y)=p(p-1)|ax+by|^{p-2}b^2
\eeq
except on the line $\{z:\,ax+by=0\}$.
Then
\beq\label{(3.3)}
\begin{array}{cl}
U_{xx}(x,y)h^2+U_{yy}(x,y)k^2
=&p(p-1)\bigg(|ax+by|^{p-2}-B_p|x|^{p-2} \bigg)h^2\\
&-p(p-1)|ax+by|^{p-2}b^2(h^2-k^2).
\end{array}
\eeq
By the property of $G(z)$, we have (see \cite{HKV})
\beq\label{(3.4)}
|ax+by|^{p-2}\leq B_p|x|^{p-2}
\eeq
in this region. So by \eqref{(3.3)} and \eqref{(3.4)} we get
\beq\label{(3.5)}
U_{xx}(x,y)h^2+U_{yy}(x,y)k^2\leq -p(p-1)|ax+by|^{p-2}b^2(h^2-k^2).
\eeq
Then \eqref{(3.2)} holds with obvious choice of $c(x,y)$.

We now consider the second case where $ z\in T$. Recall that $t_0-\varepsilon <t<t_0+\pi-\varepsilon$. We use the expression \eqref{(5.1)} for $G(z)$, then
$$
U(x,y)=r^p[\dfrac{|\cos (t_0-\theta_0)|^p}{\cos (t_0-\theta_0)}\cos(p(t-t_0)+t_0-\theta_0)-B_p \dfrac{|\cos t_0|^p}{\cos t_0}\cos(p(t-t_0)+t_0)].
$$
Since
$$
r_x=\cos t,\quad r_y=\sin t,
$$
$$
t_x=-\frac{1}{r}\sin t, \quad t_y=\frac{1}{r}\cos t,
$$
we get
$$
\begin{array}{cl}
U_{xx}(x,y)=p(p-1)r^{p-2}&\bigg(\dfrac{|\cos(t_0-\theta_0)|^p}{\cos(t_0-\theta_0)}\cos(2t-p(t-t_0)-(t_0-\theta_0))\\
&-B_p\dfrac{|\cos t_0|^p}{\cos t_0}\cos(2t-p(t-t_0)-t_0) \bigg),
\end{array}
$$
where $x=r\cos t, y=r\sin t$, $\tan \theta_0=b/a$, and
$$
U_{yy}(x,y)=-U_{xx}(x,y).
$$
Then
\beq
U_{xx}(x,y)h^2+U_{yy}(x,y)k^2=U_{xx}(x,y)(h^2-k^2).
\eeq
We claim that
\beq\label{(3.6)}
U_{xx}(x,y)\leq 0
\eeq
for $z \in T$, where $z=x+iy$.
In fact,
$$
U_{xx}(re^{it_0})=p(p-1)r^{p-2}\bigg(\dfrac{|\cos(t_0-\theta_0)|^p}{\cos(t_0-\theta_0)}\cos(t_0+\theta_0)-B_p\dfrac{|\cos t_0|^p}{\cos t_0}\cos t_0 \bigg),
$$
we know from \cite[p.249]{HKV} that
\beq\label{(3.7)}
|a\cos t_0+b\sin t_0|^{p-2}\leq B_p|\cos t_0|^{p-2},
\eeq
this is equivalent to
\beq\label{(3.8)}
|\cos (t_0-\theta_0)|^{p-2}\leq B_p|\cos t_0|^{p-2}.
\eeq
Combining \eqref{(3.8)} with the fact that
\beq\label{(3.9)}
\cos(t_0-\theta_0)\cos(t_0+\theta_0)\leq \cos ^2 t_0,
\eeq
we have
\beq\label{(3.10)}
U_{xx}(re^{it_0})\leq 0.
\eeq
Now we use the expression \eqref{(5.2)} for $G(z)$ to get
$$
U(x,y)=r^p\bigg[\dfrac{|\cos (u_0-\theta_0)|^p}{\cos (u_0-\theta_0)}\cos(p(t-u_0)+u_0-\theta_0)-B_p \dfrac{|\cos u_0|^p}{\cos u_0}\cos(p(t-u_0)+u_0)\bigg],
$$
where $u_0=t_0+\pi/p$, then
$$
\begin{array}{cl}
U_{xx}(x,y)=p(p-1)r^{p-2}&\bigg(\dfrac{|\cos(u_0-\theta_0)|^p}{\cos(u_0-\theta_0)}\cos(2t-p(t-u_0)-(u_0-\theta_0))\\
&-B_p\dfrac{|\cos u_0|^p}{\cos u_0}\cos(2t-p(t-u_0)-u_0) \bigg),
\end{array}
$$
where $x=r\cos t, y=r\sin t$, $\tan \theta_0=b/a$, so
$$
U_{xx}(re^{iu_0})=p(p-1)r^{p-2}\bigg(\dfrac{|\cos(u_0-\theta_0)|^p}{\cos(u_0-\theta_0)}\cos(u_0+\theta_0)-B_p\dfrac{|\cos u_0|^p}{\cos u_0}\cos u_0 \bigg),
$$
where $u_0=t_0+\pi/p$. We know from \cite[p.249]{HKV} that
\beq\label{(new3.7)}
|a\cos u_0+b\sin u_0|^{p-2}\leq B_p|\cos u_0|^{p-2},
\eeq
which is equivalent to
\beq\label{(3.11)}
|\cos (u_0-\theta_0)|^{p-2}\leq B_p|\cos u_0|^{p-2}.
\eeq
Combining \eqref{(3.11)} with
\beq\label{(3.12)}
\cos(u_0-\theta_0)\cos(u_0+\theta_0)\leq \cos ^2 u_0,
\eeq
we have
\beq\label{(3.13)}
U_{xx}(re^{i(t_0+\pi/p)})\leq 0.
\eeq
Write $U_{xx}(re^{it})=p(p-1)r^{p-2}u(t)$, where
$$
u(t)=A\cos |p-2|t+\textup{sgn}(2-p)B\sin |p-2|t,
$$
and
$$
A=\frac{|\cos(t_0-\theta_0)|^p}{\cos(t_0-\theta_0)}\cos\big((p-1)t_0+\theta_0\big)-B_p\frac{|\cos t_0|^p}{\cos t_0}\cos (p-1)t_0,
$$
$$
B=B_p\frac{|\cos t_0|^p}{\cos t_0}\sin (p-1)t_0-\frac{|\cos(t_0-\theta_0)|^p}{\cos(t_0-\theta_0)}\sin\big((p-1)t_0+\theta_0\big).
$$
Then $u(t)$ is a $|p-2|$-trigonometric function, thus also $|p-2|$-trigonometrically convex for $t_0<t<t_0+\pi/p$ (see \cite[p.54]{L}). We have
$U_{xx}(re^{it})=p(p-1)r^{p-2}u(t)$ is harmonic thus subharmonic within the angle $\{z=re^{it}:r>0, t_0<t<t_0+\pi/p\}$ via a direct computation.
Then, by \eqref{(3.10)}, \eqref{(3.13)} and the Phragm\'{e}n-Lindel\"{o}f theorem for subharmonic functions (see \cite{L}), we have
\beq\label{(last)}
U_{xx}(x,y)=U_{xx}(re^{it})\leq 0
\eeq
for $z \in T$. We can also use the maximum principle for harmonic functions directly to deduce  \eqref{(last)}. This  proves of \eqref{(3.6)}. Then \eqref{(3.2)} holds with $c(x,y)=-U_{xx}(x,y)$. This completes the proof of \eqref{(1.1)}.  \qed

A few comments are in order:

\medskip\noindent{\bf Remarks.}
(a)  The case $a=0$, $b=1$ of Theorem~\ref{T1.1} is contained  in  ~\cite{BW}. \\
(b)  When $a=0$, $b=1$, $p>2$, the function $U=-G$ becomes the function $U_2(x,y)$ in~\cite{BW}. \\
(c)  When $a=0$, $b=1$, $1<p<2$, the function $U=-G$ is used in  ~\cite{Pi} and in \cite{Gr}.

\section{The sharpness of the constant $B_p^{1/p}$}\label{Sharp}
To show that the constant $B_p$ is sharp, we apply a similar argument as in \cite{BW}. Let $f(z)=u(z)+iv(z)$ be analytic in the unit disc $D$ with $f(0)=0$ and $B_t$ be Brownian motion in $D$ killed upon leaving $D$. Consider the martingales $X_t=u(B_t)$ and $Y_t=v(B_t)$, we have $\langle X, Y\rangle_t=0$ and $\langle X\rangle_t-\langle Y\rangle_t=0$ (see \cite{D}). So $X$ and $Y$ are orthogonal with equal quadratic variations. Then the inequality in Theorem~\ref{T1.1} exactly reduces  to the inequality in Theorem 4.1 in \cite{HKV}.

Since $B_p^{1/p}$ is already the best constant in Theorem 4.1 of \cite{HKV}, we conclude that the constant $B_p^{1/p}$ cannot be improved in Theorem~\ref{T1.1}.

\section{Examples and applications}\label{Appl}
In this section, we give a direct application of Theorem~\ref{T1.1} to operators related to the
following discrete version of the Hilbert transform
\beq\label{(dH)}
(Da)_n= \frac{1}{\pi}\sum_{k\neq 0}\frac{a_{n-k}}{k},
\eeq
where $k$ runs over all the non-zero integers in $\mathbb Z$ and $a=(a_n)_n$. Recently, Ba\~{n}uelos and Kwa\'{s}nicki \cite{BK} proved that the operator norm of $D$ on $\ell ^p(\mathbb Z)$ is equal to the operator norm of the continuous Hilbert transform H on $L^p(\mathbb R)$. The proof in
\cite{BK}  is based on Theorem A and uses two auxilliary operators $\mathcal{J}$ [defined in \eqref{defJ}] and  $\mathcal{K}$ which satisfies $\mathcal K \mathcal J = D$. As an application of  Theorem~\ref{T1.1} we obtain the following results
concerning    $\mathcal{J}$ and  $\mathcal{K}$.

\begin{prop}\label{cor1} Let $(a_n)$ be a sequence in $\ell ^p(\mathbb Z), 1<p<\infty$. Let $\mathcal{J}a_n=\sum\limits_{m\in \mathbb Z}\mathcal{J}_ma_{n-m}$, where
\begin{equation}\label{defJ}
\mathcal{J}_n=\frac{1}{\pi n}\bigg(1+ \int_0^\infty \frac{2y^3}{(y^2+\pi^2n^2)\sinh^2y }dy \bigg)
\end{equation}
for $n\neq 0$, and $\mathcal{J}_0=0$. Then for $a, b \in \mathbb R$,
\beq\label{(cor1.1)}
\|(aI+b\mathcal{J})a_n\|_p\leq B_p^{1/p}\|a_n\|_p,
\eeq
where $B_p$ is given by \eqref{(1.0)} and $I$ is the identity operator: the convolution with kernel $I_0=1, I_n=0$ for $n\neq 0$. The constant  $B_p^{1/p}$ is the best possible in this inequality.
\end{prop}

\emph{Proof.} We use the notation in \cite{BK}. We only need to redefine the operator in (2.5) in \cite{BK} that
\beq\label{(corproof1.1)}
J_Aa_n=\mathbb E_{(x_0,y_0)}\big(a\|A\|M_{\zeta -} + bA\star M_{\zeta -} | Z_{\zeta -}=(2\pi n, 0)\big).
\eeq
Since the conditional expectation is a contraction on $L^p, 1<p<\infty$, it follows from \eqref{(1.1)} in Theorem~\ref{T1.1} that
\beq\label{(corproof1.2)}
\|J_Aa_n\|_p\leq B_p^{1/p}\|A\|\|a_n\|_p.
\eeq
Let
$$
H=\left[ \begin{array}{ccc} 0 & -1 \\ 1 & 0 \end{array} \right],
$$
we have
\beq\label{(corproof1.3)}
\|J_Ha_n\|_p\leq B_p^{1/p}\|a_n\|_p.
\eeq
Notice that $\mathbb E_{(x_0,y_0)}\big(M_{\zeta -} | Z_{\zeta -}=(2\pi n, 0)\big)=Ia_n$, where $I$ is the identity operator, then following the same proof in \cite{BK}, we deduce \eqref{(cor1.1)}.

The sharpness of the constant is due to the sharpness of Proposition~\ref{cor2} and the fact that
\beq\label{(cor1.4)}
\|(a \mathcal{K}+b D)a_n\|_p\leq \|(a I + b \mathcal{J} )a_n\|_p
\eeq
for any sequence $(a_n)\in \ell ^p(\mathbb Z), 1<p<\infty$ and $a, b \in \mathbb R$.
\qed

\begin{prop}\label{cor2} Let $(a_n)$ be a sequence in $\ell ^p(\mathbb Z), 1<p<\infty$, $D$ be defined in \eqref{(dH)}. Let $\mathcal{K}$ be the convolution operator in Section 2.3 in \cite{BK} with kernel $(\mathcal{K}_n)$ such that $\mathcal{K}_n\geq 0$ for all $n$ and the sum of all $\mathcal{K}_n$ is equal to 1. Then for $a, b \in \mathbb R$,
\beq\label{(cor1.2)}
\|(a \mathcal{K}+b D)a_n\|_p\leq B_p^{1/p}\|a_n\|_p,
\eeq
where $B_p$ is given by \eqref{(1.0)}. The constant $B_p^{1/p}$ is the best possible in this inequality.
\end{prop}

\emph{Proof.} By Section 2.3 in \cite{BK},
$$
Da_n=\mathcal{K}\mathcal{J}a_n,
$$
then by Proposition~\ref{cor1}, we have
\beq\label{(cor2.1)}
\|(a \mathcal{K}+b D)a_n\|_p=\|\mathcal{K} (a I + b \mathcal{J} )a_n\|_p
\leq B_p^{1/p}\|a_n\|_p.
\eeq

To deduce the sharpness, we define the dilation operators $T_\varepsilon$ for any $\varepsilon >0$ and $1<p<\infty$ by $(T_\varepsilon f)(x)=\varepsilon^{1/p}f(\varepsilon x)$, then $\|T_\varepsilon\|_{p,p}=1$ for all $\varepsilon >0$. Notice that $\mathcal{K}$ is a convolution operator with kernel $(\mathcal{K}_n)$ such that $\mathcal{K}_n\geq 0$ for all $n$ and $\sum_{n\in \mathbb Z}\mathcal{K}_n=1$ (see \cite{BK}).  Because of Theorem 4.2 in \cite{La}, we can work on the real line and replace $D$ and $\mathcal{K}$ by
\beq\label{(cor2.2)}
(M_D f)(x)= \mathrm{p.v.}\frac{1}{\pi}\sum_{m\neq 0}\frac{f(x-m)}{m}
\eeq
and
\beq\label{(cor2.3)}
(M_{\mathcal{K}} f)(x)=\sum_{m\in \mathbb Z}\mathcal{K}_m f(x-m),
\eeq
respectively. It is known by \cite{La} that
$$
\lim_{\varepsilon\rightarrow 0}(T_{1/\varepsilon}M_D T_\varepsilon f)(x)=(Hf)(x),
$$
where $H$ is the Hilbert transform. We claim that
\beq\label{(cor2.4)}
\lim_{\varepsilon\rightarrow 0}(T_{1/\varepsilon}M_\mathcal{K} T_\varepsilon f)(x)=(If)(x),
\eeq
for $\mathrm{a.e.}\; x\in \mathbb R$ and $f\in L^p(\mathbb R)$, where $I$ is the identity operator such that $(If)(x)=f(x)$. In fact, for any $f\in \mathcal{S}(\mathbb R)$ (Schwartz function), we have
$$
\begin{array}{cl}
&\lim\limits_{\varepsilon\rightarrow 0}(T_{1/\varepsilon}M_\mathcal{K} T_\varepsilon f)(x)\\
=&\lim\limits_{\varepsilon\rightarrow 0}\sum_{|m|\leq N}\mathcal{K}_m f(x-\varepsilon m)
+\lim\limits_{\varepsilon\rightarrow 0}\sum_{|m|> N}\mathcal{K}_m f(x-\varepsilon m)\\
=&\sum_{|m|\leq N}\mathcal{K}_m f(x)+\lim\limits_{\varepsilon\rightarrow 0}\sum_{|m|> N}\mathcal{K}_m f(x-\varepsilon m)
\end{array}
$$
for any $N>0$. Then
$$
\begin{array}{cl}
&\big|f(x)-\lim\limits_{\varepsilon\rightarrow 0}(T_{1/\varepsilon}M_\mathcal{K} T_\varepsilon f)(x) \big|\\
=&\big|\sum_{m\in \mathbb Z}\mathcal{K}_m f(x)-\sum_{|m|\leq N}\mathcal{K}_m f(x)-\lim\limits_{\varepsilon\rightarrow 0}\sum_{|m|> N}\mathcal{K}_m f(x-\varepsilon m)\big|\\
=&\big|\sum_{|m|> N}\mathcal{K}_m f(x)-\lim\limits_{\varepsilon\rightarrow 0}\sum_{|m|> N}\mathcal{K}_m f(x-\varepsilon m) \big|\\
\leq &\lim\limits_{\varepsilon\rightarrow 0}\sum_{|m|> N}\mathcal{K}_m\big|f(x)-f(x-\varepsilon m)\big|\\
\leq & C(f)\sum_{|m|> N}\mathcal{K}_m
\end{array}
$$
for any $N>0$. Since $\mathcal{K}_n\geq 0$ for all $n$ and $\sum_{n\in \mathbb Z}\mathcal{K}_n=1$, letting $N\rightarrow \infty$, we get
$$
\lim\limits_{\varepsilon\rightarrow 0}(T_{1/\varepsilon}M_\mathcal{K} T_\varepsilon f)(x)=f(x)
$$
for $x\in \mathbb R$, $f\in \mathcal{S}(\mathbb R)$. For $f\in L^p(\mathbb R)$, we can get \eqref{(cor2.4)} for $\mathrm{a.e.}\; x\in \mathbb R$ by the standard density argument.

Then, we have
$$
\|aI+bH\|_{p,p}\leq \sup_{\varepsilon}\|T_{1/\varepsilon}(aM_{\mathcal{K}}+bM_D)T_\varepsilon\|_{p,p}\leq \|aM_{\mathcal{K}}+bM_D\|_{p,p}.
$$
By Theorem 4.2 in \cite{La},
$$
\|aM_{\mathcal{K}}+bM_D\|_{p,p}=\|a\mathcal{K} + b D\|_{p,p},
$$
so we have
$$
\|a\mathcal{K} + b D\|_{p,p}=\|aI+bH\|_{p,p}=B_p^{1/p}.
$$
This finishes the proof of Proposition~\ref{cor2}.

\qed

For the operator $aI+bD$, $a, b\in \mathbb R$, applying the method in \cite[Lemma 4.3]{La}, we immediately obtain
$$
\|aI+bD\|_{p,p}\geq \|aI+bH\|_{p,p}=B_p^{1/p}.
$$

We conjecture that $\|aI+bD\|_{p,p}=\|a \mathcal{K}+b D\|_{p,p}=B_p^{1/p}$. The
solution of this conjecture may require   additional ideas as $I$ and $D$ are natural projections
of nonorthogonal  martingales.

\end{document}